\documentclass[aps,preprint,showpacs,superscriptaddress,groupedaddress]{revtex4-1}  
\usepackage{graphicx}  
\usepackage{dcolumn}   
\usepackage{bm}        
\usepackage{amssymb}   
\usepackage{amsmath}   
\usepackage{listings}
\usepackage{array}
\usepackage{graphicx}
\usepackage{lineno}
\usepackage{dcolumn}
\usepackage{color}
\usepackage{overpic}
\usepackage{multirow}
\usepackage{hyperref}

\begin{document}

\widetext


\title{Enhanced inequalities about arithmetic and geometric means}
\author{Fang Dai$^{1}$, Li-Gang Xia \\ 1. Futian middle school, Ji'an City, Jiangxi Province, China, 343000}

\begin{abstract}
    For $n$ positive numbers ($a_k$, $1\leq k \leq n$), enhanced inequalities about the arithmetic mean ($A_n \equiv \frac{\sum_ka_k}{n}$) and the geometric mean ($G_n\equiv \sqrt[n]{\Pi_ka_k}$) are found if some numbers are known, namely,
\begin{equation}
    \frac{G_n}{A_n} \leq (\frac{n-\sum_{k=1}^mr_k}{n-m})^{1-\frac{m}{n}}(\Pi_{k=1}^mr_k)^{\frac{1}{n}}\leq 1 \:, \nonumber
\end{equation}
if we know $a_k=A_nr_k$ ($1\leq k\leq m\leq n$) for instance, and
\begin{equation}
    \frac{G_n}{A_n} \leq \frac{1}{(1-\frac{m}{n})\Pi_{k=1}^mr_k^{\frac{-1}{n-m}}+\frac{1}{n}\sum_{k=1}^mr_k}\leq 1 \: ,\nonumber 
\end{equation}
if we know $a_k=G_nr_k$ ($1\leq k\leq m \leq n$) for instance. These bounds are better than those derived from S.~H.~Tung's work~\cite{Tung}.
\end{abstract}
\maketitle 

\section{Introduction}
Let $a_1,a_2,\ldots,a_n$ denote $n$ positive numbers. Let $A_n$ be their arithmetic mean, $\frac{\sum_ka_k}{n}$, and $G_n$ be their geometric mean, $\sqrt[n]{\Pi_ka_k}$. We shall prove the following inequalities.
\begin{equation}\label{eq:xia1}
    \frac{G_n}{A_n} \leq (\frac{n-\sum_{k=1}^mr_k}{n-m})^{1-\frac{m}{n}}(\Pi_{k=1}^mr_k)^{\frac{1}{n}}\leq 1 \: .
\end{equation}
if we know $a_k=A_nr_k$ ($1\leq k\leq m\leq n$) for instance.
\begin{equation}\label{eq:xia2}
    \frac{G_n}{A_n} \leq \frac{1}{(1-\frac{m}{n})\Pi_{k=1}^mr_k^{\frac{-1}{n-m}}+\frac{1}{n}\sum_{k=1}^mr_k} \leq 1 
\end{equation}
if we know $a_k=G_nr_k$ ($1\leq k\leq m \leq n$) for instance.

S. H. Tung~\cite{Tung} obtained the following lower bound for $A_n-G_n$ in terms of the smallest value, $a$, and the largest value, $A$.
\begin{equation}\label{eq:Tung0}
    A_n-G_n \geq \frac{(\sqrt{A}-\sqrt{a})^2}{n} \: . 
\end{equation}
We will see that our results are better than this bound.

\section{Proof of the inequalities}
Suppose we know the first $m$ numbers, $a_1,a_2,\ldots,a_m$. We can construct the following inequality. 
\begin{eqnarray}
    && \sum_{k=1}^m\lambda_k^{n-m}a_k + \frac{1}{\Pi_{i=1}^m\lambda_i}\sum_{k=m+1}^n a_k \\
    = && \sum_{k=1}^m(\lambda_k^{n-m}-\frac{1}{\Pi_{i=1}^m\lambda_i})a_k + \frac{1}{\Pi_{i=1}^m\lambda_i}\sum_{k=1}^n a_k \label{eq:step2}\\
    \geq && n\sqrt[n]{\Pi_ka_k}
\end{eqnarray}
Suppose we know $a_k=A_nr_k$ ($k=1,2,\ldots,m$). Inserting them into Eq.~\ref{eq:step2}, we obtain an upper bound of $G_n/A_n$ as a function of $\lambda_1,\lambda_2,\ldots,\lambda_m$.
\begin{equation}
    \frac{G_n}{A_n} \leq \frac{1}{n}\sum_{k=1}^m\lambda_k^{n-m}r_k+\frac{1}{\Pi_{k=1}^m\lambda_k}\frac{n-\sum_{k=1}^mr_k}{n} \equiv f(\lambda_1,\lambda_2,\ldots,\lambda_m) \: . 
\end{equation}
$\frac{\partial f}{\partial \lambda_i}=0$ ($i=1,2,\ldots,m$) gives the best choice, namely,
\begin{eqnarray}
    && \lambda_i = \left(\frac{n-\sum_{k=1}^mr_k}{n-m}\Pi_{k=1}^mr_k^{\frac{1}{n-m}}\right)^{\frac{1}{n}}r_i^{\frac{-1}{n-m}} \: , \: (i=1,2,\ldots,m)
\end{eqnarray}
and hence the bound in Ineq.~\ref{eq:xia1}.

Similarly suppose we know $a_k=G_nr_k$ ($k=1,2,\ldots,m$). Inserting them into Eq.~\ref{eq:step2}, we obtain an upper bound of $G_n/A_n$ as a function of $\lambda_1,\lambda_2,\ldots,\lambda_m$.
\begin{equation}
    \frac{G_n}{A_n} \leq \frac{1}{\Pi_{k=1}^m\lambda_k(1-\frac{1}{n}\sum_{k=1}^mr_k\lambda_k^{n-m})+\frac{1}{n}\sum_{k=1}^mr_k} \equiv g(\lambda_1,\lambda_2,\ldots,\lambda_m) \: . 
\end{equation}
$\frac{\partial f}{\partial \lambda_i}=0$ ($i=1,2,\ldots,m$) gives the best choice, namely,
\begin{eqnarray}
    && \lambda_i = r_i^{\frac{-1}{n-m}} \: , \: (i=1,2,\ldots,m)
\end{eqnarray}
and hence the bound in Ineq.~\ref{eq:xia2}.

\begin{figure}
    \includegraphics[width=0.45\textwidth]{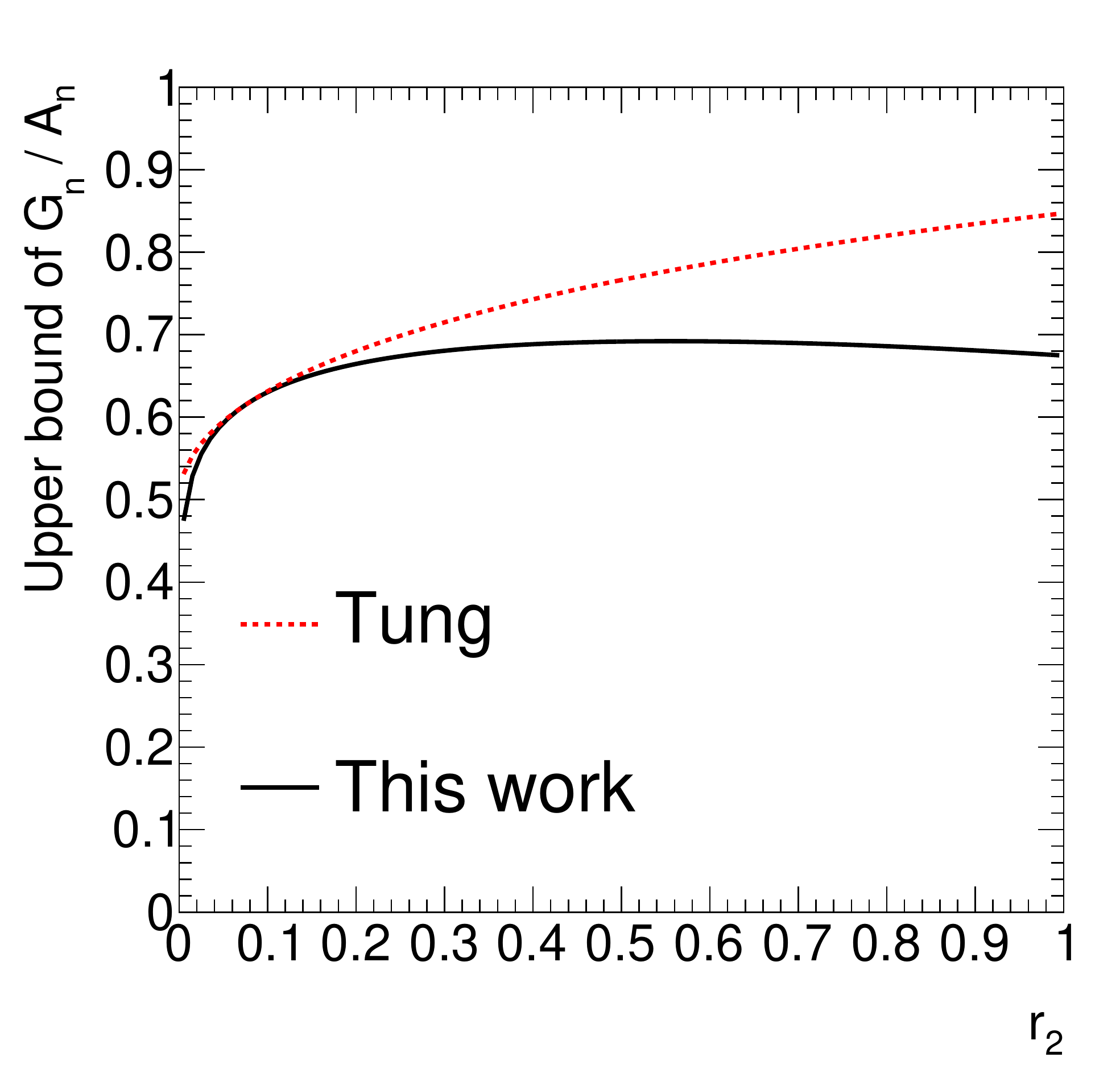}
    \includegraphics[width=0.45\textwidth]{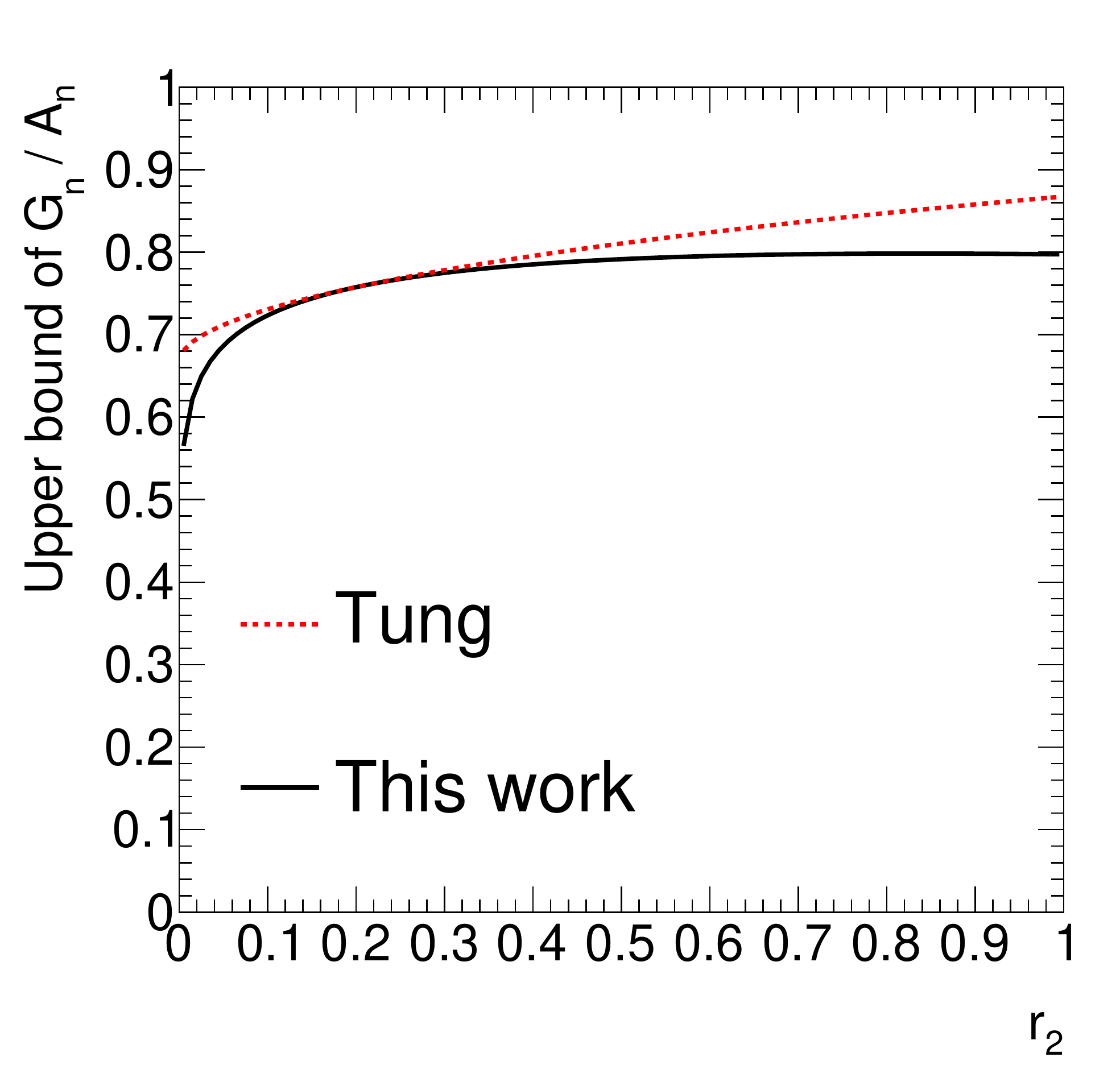}
    \caption{\label{fig:comparison}
        Comparison of different upper bounds of $G_n/A_n$ as a function of $r_2$ for $r_1=5$ and $n=5$. The left diagram correponds to the case of $A=A_nr_1$ and $a=A_nr_2$ while the right diagram corresponds to the case of $A=G_nr_1$ and $a=G_nr_2$. The results derived from Ref.~\cite{Tung} are represented by the red-dashed curves (namely, Ineq.~\ref{eq:tung1} and ~\ref{eq:tung2}. The black curves represent our results
        (namely, Ineq.~\ref{eq:xia1} and ~\ref{eq:xia2}).
    }
\end{figure}

For comparison, Tung's inequality can be written into the following form,
\begin{equation}
    \frac{G_n}{A_n} \leq 1 - \frac{1}{n}(\sqrt{r_1}-\sqrt{r_2})^2 \:, \label{eq:tung1}
\end{equation}
if we know $A=A_nr_1$ and $a=A_nr_2$ with $0<r_2\leq 1\leq r_1 \leq n-r_2$, or
\begin{equation}
    \frac{G_n}{A_n} \leq \frac{1}{1+\frac{1}{n}(\sqrt{r_1}-\sqrt{r_2})^2} \:, \label{eq:tung2}
\end{equation}
if we know $A=G_nr_1$ and $a=G_nr_2$ with $0<r_2\leq 1\leq r_1$. We can show that our results are better than these bounds using simple calculus. For illustration, different upper bounds of $G_n/A_n$ as a function of $r_2$ with $r_1=5$ and $n=10$ are compared in Fig.~\ref{fig:comparison}.

\end{document}